\pdfoutput=1

\documentclass{amsart}

\usepackage{snapshot}

\usepackage{amsmath}
\usepackage{mathrsfs} 	% for \mathscr
\usepackage{amssymb}
\usepackage{amsthm} 	% provides Theorem environment
\usepackage{thmtools} 	% for manual numbering of theorems
\usepackage{amscd} 		% provides CD environment for commutative diagrams (and exact seqs)

\usepackage{subcaption}
\usepackage{graphicx}
\DeclareGraphicsExtensions{.pdf,.eps,.png,.jpg}
\usepackage{cite}
\usepackage{hyperref}

\newtheorem{lemma}{Lemma}
\newtheorem{thm}{Theorem}
\newtheorem{prop}[thm]{Proposition}
\theoremstyle{remark} \newtheorem{remark}{Remark}
\theoremstyle{definition} \newtheorem{ex}{Example} \newtheorem*{defn}{Definition}

\newcommand{\bbA}{\mathbb{A}}

\newcommand{\bbN}{\mathbb{N}}

\newcommand{\calO}{\mathcal{O}}

\newcommand{\scrI}{\mathscr{I}}

\newcommand{\wt}[1]{\widetilde{#1}} 
\newcommand{\setrow}[2]{\{#1_1, \dots, #1_{#2}\}}
\newcommand{\row}[2]{#1_1, \dots, #1_{#2}}

\newcommand{\Spec}{\operatorname{Spec}}
\newcommand{\Proj}{\operatorname{Proj}}
\newcommand{\codim}{\operatorname{codim}}

\title{Monomial principalization in the singular setting}

\author{Corey Harris}
\address{Department of Mathematics, Florida State University, Tallahassee FL, 32306, USA}
\email{charris@math.fsu.edu}
\urladdr{\url{http://www.math.fsu.edu/~charris}}
\thanks{Thanks to Paolo Aluffi for suggesting this problem and for his patience and feedback while this work was prepared.}

\begin{document}

\begin{abstract}
	We generalize an algorithm by Goward for principalization of monomial ideals in nonsingular varieties to work on any scheme of finite type over a field. The normal crossings condition considered by Goward is weakened to the condition that components of the generating divisors meet as complete intersections.
	This leads to a substantial generalization of the notion of monomial scheme; we call the resulting schemes {\em `c.i.~monomial'.\/} We prove that c.i.~monomial schemes in arbitrarily singular varieties can be principalized by a sequence of blow-ups at codimension~2 c.i.~monomial centers.
\end{abstract}

\maketitle

\section{Introduction}
\label{sec:intro}
	{\em Monomial schemes\/} are schemes defined as intersections of collections 
	of components from a fixed normal crossing divisor in a nonsingular variety. 
	In [Gow05], Goward proves that monomial schemes may be principalized by a 
	sequence of blow-ups along codimension~2 monomial schemes. In the singular 
	setting, this definition of monomial scheme is not available, because the 
	notion of normal crossings requires a nonsingular ambient variety. We 
	consider a much weaker condition, which makes no nonsingularity assumption on 
	the ambient variety: essentially, divisors meet with `c.i.~crossings' if they 
	intersect along subschemes with the expected dimension. (For example, any two 
	distinct irreducible curves on a smooth surface meet with c.i.~crossings, 
	regardless of whether they are nonsingular or meet transversally.) This leads 
	to a generalization of the notion of monomial scheme, {`c.i.~monomial 
	schemes'.\/} See \autoref{sec:defs_and_examples} for formal definitions.

	We extend Goward's result to the c.i.~monomial case, showing that this much larger class of subschemes can be principalized via Goward's procedure.

	Our result has been used in recent work on computations of Segre classes, cf.~\cite{Alu}, Theorem~1.1.

\section{Definitions and Examples}
\label{sec:defs_and_examples}

	Throughout, $X$ will denote a scheme of finite type over an arbitrary field.
	By regular sequence, we mean a sequence $x_1,\dots,x_n$ of elements in a ring 
	$R$ such that $(x_1,\dots,x_n) \subset R$ is a proper ideal and, for each 
	$i$, the image of $x_i$ in $R/(\row x{i-1})$ is a non-zerodivisor, see 
	\cite[p. 243]{Eis95}.

	\begin{defn}[c.i.~crossings]
		Let $\row Yn \subset X$ be Cartier divisors. We say that $\setrow Yn$ has
		\emph{c.i.~crossings} if for every subset $A \subset \setrow Yn$ and 
		every point $p \in \cap_{Y \in A} Y$, the local equations $y_i$ for the 
		$Y \in A$ form a regular sequence at $p$.
	\end{defn}

	Note that the definition requires each $Y_i$ to be cut out locally by a 
	non-zerodivisor, making $Y_i$ an effective Cartier divisor in $X$.  
	Note also that the condition places no restrictions on $X$.

	The following definition will be used only in the introduction to compare concepts.

	\begin{defn}[simple normal crossings]
		Let $\row Yn \subset X$ be Cartier divisors.  We say that $\setrow Yn$ has \emph{simple normal crossings} if for every subset $A \subset \setrow Yn$, the intersection $Z = \cap_{Y \in A} Y$ is nonsingular with $\codim_X Z = \# A$. If $D = \sum a_i Y_i$, with $a_i \geq 0$, we say $D$ is
		a \emph{simple normal crossings divisor} or \emph{s.n.c~divisor}.
	\end{defn}

	The s.n.c.~condition on singletons requires each $Y_i$ to be nonsingular, and the condition on the empty set means $X$ itself must be nonsingular.

\begin{figure}[hbt]
	 	\begin{subfigure}{.49\linewidth}
			\centering
			\includegraphics[width=.97\linewidth]{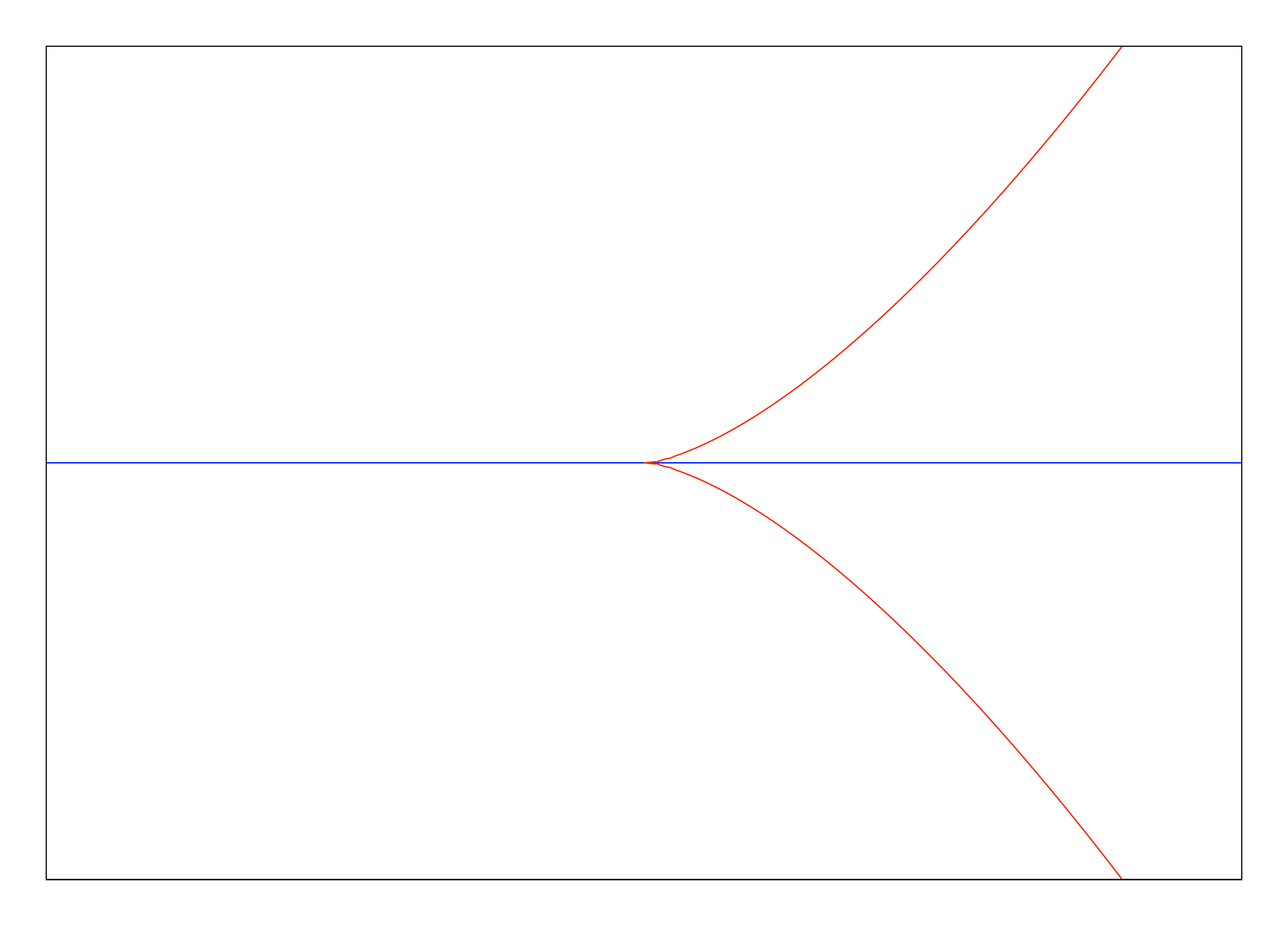}
			\caption{Curves $y$, $y^2-x^3$ in the plane}
			\label{fig:cuspidalCubicWithAxis}
		\end{subfigure}
		~
		\begin{subfigure}{.49\linewidth}
			\centering
			\includegraphics[width=.9\linewidth]{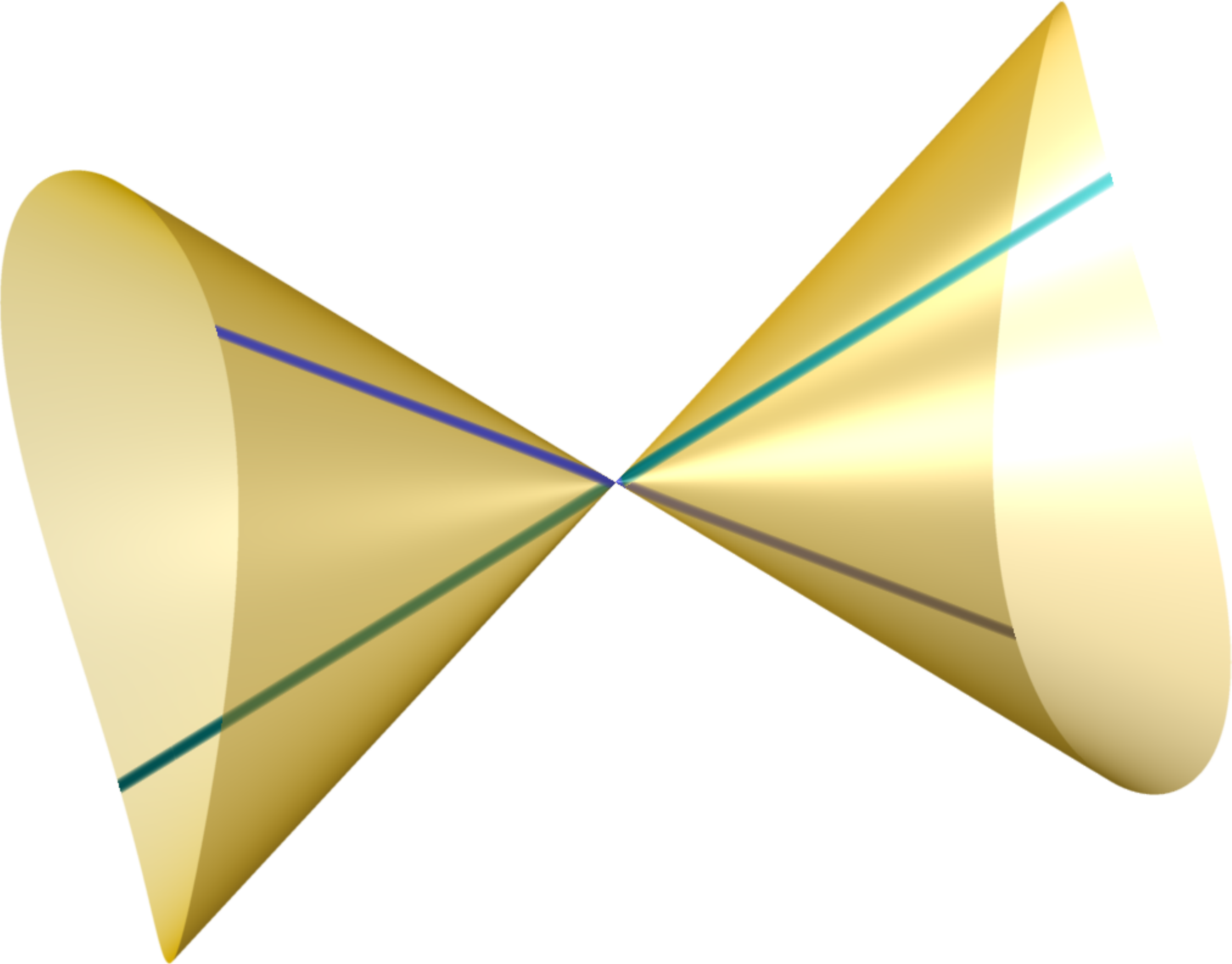}
			\caption{Curves $x$, $y$ on $xy-z^2$}
			\label{fig:coneWithAxes}
		\end{subfigure}
		\caption{Examples}
	\end{figure}

	\begin{ex}
		Let $\bbA_k^2 = \Spec k[x,y]$ and consider the curves $Y_1$ defined by $y$ and $Y_2$ defined by $y^2-x^3$. These two curves do not meet with simple normal crossings because $y^2-x^3$ is not smooth in the intersection.  So $Y_1 + Y_2$ is not an s.n.c. divisor, but $\{Y_1, Y_2\}$ does have c.i.~crossings since $y^2-x^3$ is sent to a non-zerodivisor in the integral domain $k[x,y]/(y)$. Observe (\autoref{fig:cuspidalCubicWithAxis}) that these curves do not meet transversally.
	\end{ex}

	\begin{ex}
		Now consider the cone $C = \Spec k[x,y,z]/(xy-z^2)$ and the subschemes $Y_x,Y_y$ cut out by the ideals $(x)$ and $(y)$, respectively.  Since $C$ is singular, there are no s.n.c. divisors in sight.  However $(k[x,y,z]/(xy-z^2))/(y) \cong k[x,z]/(z^2)$, so $y, x$ forms a regular sequence at every point in $C$.  Thus $\{Y_x, Y_y\}$ has c.i.~crossings in $C$.
	\end{ex}

	\begin{remark}
	\label{rmk:localization}
		In both of the above examples, we checked the c.i.~crossing condition only affine-locally, even though the definition is in terms of stalks.  This is sufficient because if $\row xn$ is a regular sequence in a Noetherian ring $R$, then it remains a regular sequence in $R_p$ for any prime ideal $p \subset R$ containing $\setrow xn$.
	\end{remark}

	\begin{defn}
		A subscheme $Z \subset X$ is called \emph{monomial} if it is cut out by effective divisors which are supported on a fixed s.n.c.~divisor.  By analogy, if $\setrow Yn$ has c.i.~crossings, we will say $Z$ is a \emph{c.i.~monomial subscheme} (with respect to $\setrow Yn$) if $Z$ is cut out by divisors of the form $\sum a_i Y_i$ with $a_i \geq 0$. As well, if $\beta: \wt X \rightarrow X$ is the blowup of $X$ at a c.i.~monomial subscheme, we will call $\beta$ (or just $\wt X$) a \emph{c.i.~monomial blowup} (with respect to $\setrow Yn$).
	\end{defn}

	We caution the reader that c.i.~monomial subschemes are \emph{not} 
	necessarily complete intersections (nor monomial).

	\newcommand{\stateMainThm}{ % we make the statement a macro, so it can be called again later
		Let $\setrow Yn$ have c.i.~crossings. If $D_1, \dots, D_h$ are given by 
		$D_j = \sum a_{ij} Y_i$, where $a_{ij} > 0$, then there exists a sequence of c.i.~monomial blowups at codimension 2 centers
		\[
			\begin{CD}
				\wt X = X_n @>\beta_{n}>> X_{n-1} @>\beta_{n-1}>> \cdots 
				@>\beta_{2}>> X_1 @>\beta_{1}>> X
			\end{CD}
		\]
		such that $(\scrI_{D_1} + \dots + \scrI_{D_h}) \calO_{X_i}$ is c.i.~monomial for each $i$ and $(\scrI_{D_1} + \dots + \scrI_{D_h}) \calO_{\wt X}$ is locally principal.
	}
	\begin{thm}[Main Theorem]\label{mainThm}
		\stateMainThm
	\end{thm}

	Goward's theorem \cite[Theorem 2]{Gow05} is the analogous statement to \autoref{mainThm} for monomial subschemes.  The algorithm here is a direct generalization of Goward's, and our method of proof follows his.  We first need to know that c.i.~monomial blowups preserve c.i.~crossings and so we verify this in \autoref{sec:blowingUp}.  Following this, we give the proof of \autoref{mainThm} in \autoref{sec:mainThmProof}.

\section{C.i. monomial blowups preserve c.i.~crossings}\label{sec:blowingUp}

	Assume $\setrow Yn$ has c.i.~crossings and $n \geq 2$.  Let 
	$D_1 = \sum a_i Y_i$ and $D_2 = \sum b_i Y_i$ with $a_i, b_j \geq 0$. If 
	${\beta: \wt X \rightarrow X}$ is the blowup of $X$ at $Y_1 \cap Y_2$, let 
	$E$ denote the exceptional divisor in $\wt X$ and let $\wt Y_i$ denote the 
	proper transform of $Y_i$ in $\wt X$. We will verify that $\beta^* D_i$, 
	which is supported on $\{E, \wt{Y}_1, \dots, \wt{Y}_n\}$, has c.i.~crossings.

	The proof will be by induction on the number of divisors.  The base case is 
	handled in \autoref{blowupBaseCaseProp} and the inductive step handled in 
	\autoref{blowupProp}.

	\begin{prop}\label{blowupBaseCaseProp}
		Suppose $\{Y_1, Y_2\}$ has c.i.~crossings on $X$.  Let 
		$\beta: \wt X \rightarrow X$ denote the blowup of $X$ at $Y_1 \cap Y_2$. 
		Then $\{E, \wt Y_1, \wt Y_2\}$ has c.i.~crossings in $\wt X$. 

		\begin{proof}
			There are eight subsets of $\{\wt{Y_1}, \wt{Y_2}, E\}$.  The 
			intersections corresponding to 
			$\emptyset$, $\{\wt{Y_1}, \wt{Y_2} \}$, $\{ \wt{Y_1}, \wt{Y_2}, E \}$ 
			are empty, and those corresponding to $\{\wt{Y_1}\}$, $\{\wt{Y_2}\}$, 
			$\{E\}$ are effective Cartier. This leaves the intersections 
			corresponding to $\{\wt{Y_1}, E\}$ and $\{\wt{Y_2}, E\}$. We show the 
			result for $\{\wt{Y_1}, E\}$.

			It suffices to prove the result affine-locally. Since $Y_1$, $Y_2$ 
			are effective Cartier divisors on $X$, there is an affine open cover 
			$\{ U_\alpha = \Spec R_\alpha\}$ of $X$ such that $y_{1,\alpha}$ is 
			a local equation for $Y_1$ on $U_\alpha$ and $y_{2,\alpha}$ is a 
			local equation for $Y_2$. Now let $U = \Spec R$ be a member of such 
			a cover and let $(y_1),(y_2) \subset R$ be principal ideals defining 
			$Y_1$ and $Y_2$, respectively.  Then the blowup $\wt U$ of $U$ 
			centered at $Y_1 \cap Y_2$ is $\Proj R[a_1,a_2]/(a_2 y_1 - a_1 y_2)$, 
			see \cite[B.6.10]{Ful98}. Consider the open set 
			${D(a_2) := \{[p] \in \Spec R \;|\; a_2 \notin p\}}$ which we can 
			write as
			\[
				D(a_2) = \Spec R[a_1]/(y_1 - a_1 y_2).
			\]  
			The pullback $\beta^* Y_1$ is cut out by $a_1 y_2$ in $D(a_2)$ and the exceptional divisor is cut out by $y_2$, so the proper transform of $Y_1$ is cut out in $D(a_2)$ by $a_1$.  Since ${(R[a_1]/(y_1-a_1 y_2))/(a_1) \cong R/(y_1)}$, and $y_2$ is not a zerodivisor in $R/(y_1)$ by assumption, we have that $a_1, y_2$ is a regular sequence in $R[a_1]/(y_1 - a_1 y_2)$ corresponding to the intersection $Y_1 \cap Y_2$.

			The proof is completed by noting again (\autoref{rmk:localization}) that localization preserves regular sequences, so that $a_1, y_2$ is a regular sequence in $\calO_{\wt X, p}$ for each $p \in \wt {Y}_1 \cap E$.
		\end{proof}
	\end{prop}

	We note that as a scheme of finite type over a Noetherian ring, $X$ is Noetherian. In particular $\calO_{X,p}$ is a Noetherian local ring for all $p \in X$. In the proof of \autoref{blowupProp}, we will need the following result.

	\begin{lemma}\label{regSeqLemma}
		If $R$ is a Noetherian local ring and $x_1, \dots, x_r$ is a regular sequence of elements	in the maximal ideal of $R$, then any permutation of $x_1, \dots, x_r$ is again a regular sequence.
		\begin{proof}
			See \cite[Theorem 17.2]{Eis95}.
		\end{proof}
	\end{lemma}

	We will make repeated use of the lemma, along with the idea that if $\row xn$ is a regular sequence, then $\row xk$ is a regular sequence for $1 \leq k \leq n$.

	\begin{prop}\label{blowupProp}
		Suppose $\setrow Yn$ has c.i.~crossings on $X$.  Let 
		$\beta: \wt X \rightarrow X$ denote the blowup of $X$ at $Y_1 \cap Y_2$. 
	 	Then $\{E, \wt{Y}_1, \dots, \wt{Y}_n\}$ has c.i.~crossings in $\wt X$. 
	 	\begin{proof}
	 		By induction, assume $\{E, \wt{Y}_1, \dots, \wt{Y}_{n-1}\}$ has c.i.~crossings. We will show that if $A$ contains any (possibly empty) subset of $\{E,\wt Y_1, \wt Y_2\}$, then the corresponding intersection is cut out	by a regular sequence at each point in an affine chart on the blowup.

	 		Let $A \subset \{E, \row {\wt Y}n\}$.  We want $\cap_{Z \in A} Z \neq \emptyset$ so $\{\wt Y_1, \wt Y_2\}$ cannot be a subset of $A$. Assume then without loss of generality that $\wt Y_2 \notin A$.  Assume also that $\wt Y_n \in A$. Each element of $\setrow Yn$ is an effective Cartier divisor so again we can find a cover of $X$ by open affines $\{U_\alpha = \Spec R_\alpha\}$ such that for each $i$, and each $Y \in \setrow Yn$, we have a non-zerodivisor $y \in R_\alpha$ which cuts out $Y$ in $U_\alpha$.  Choose such a $U = \Spec R$ and let $y_i \in R$ be a local equation for $Y_i$.  We can then write the blowup of this chart $\wt U = \Proj R[a_1,a_2]/(a_2 y_1 - a_1 y_2)$. As before, we work affine-locally in $D(a_2) = \Spec R[a_1]/(y_1 - a_1 y_2)$, and we have that $E$ is cut out by $y_2$ and $\wt Y_1$ is cut out by $a_1$.

	 		Assume that $A$ contains both $E$ and $\wt Y_1$.  By induction we can find a regular sequence for $A \backslash \{\wt Y_n\}$ in the elements $\{y_2, a_1, \row sk\}$ (where ${\# (A \backslash \{\wt Y_n\}) = k+2}$) and so \autoref{regSeqLemma} ensures that $a_1, y_2, \row sk$ is a regular sequence at each \newline ${ p \in \cap_{Y \in A\backslash \{\wt Y_n\}} Y }$. Let $r$ be a local equation for $\wt Y_n$ on $D(a_2)$ such that we have a regular sequence in the elements $r, y_1, y_2, \row sk$ at each point $p \in \cap_{Y \in A} Y$. Then by \autoref{regSeqLemma} these elements form a regular sequence in any order.  Thus $y_2, r, \row sn$ form a regular sequence in $R/(t_1)$, showing that $a_1, y_2, r, \row sn$ is a regular sequence at each point $p \in \cap_{Y \in A} Y$.

	 		If $A$ contains $E$ but does not contain $\wt Y_1$, then we can use \autoref{regSeqLemma} again to get the regular sequence $\row sn, r,y_2,a_1$ which shows that $\row sn, r, y_2$ is a regular sequence.  Similarly, if $A$ contains $\wt Y_1$ but not $E$, we can rearrange to get $\row sn, r,a_1,y_2$, which shows $\row sn, r,a_1$ is a regular sequence.  Truncating this sequence again shows that $\row sn, r$ is a regular sequence.  This is the case where
	 		$A$ contains neither $E$ nor $\wt Y_1$.
	 	\end{proof}
	\end{prop}

	In this paper, we will only need the result as stated in \autoref{blowupProp}, but since the fact is true in greater generality, we provide \autoref{blowupPropforallCodim} for completeness.  In the following proof we use \autoref{blowupProp} as the base case and induct on the number of divisors that cut out the center for the blowup.

	\begin{prop}\label{blowupPropforallCodim}
		Suppose $\setrow Yn$ has c.i.~crossings on $X$.  Let ${\beta: \wt X \rightarrow X}$ denote the blowup of $X$ at $\cap_{i=1}^r Y_i$, where $r \leq n$. Then $\{E, \row {\wt Y}n\}$ has c.i.~crossings in $\wt X$.
	 	\begin{proof}
	 		Let $A \subset \{E, \row {\wt Y}n\}$ with $\wt {Y}_1 \in A$.  We want $\cap_{Z \in A} Z \neq \emptyset$ so we require $A \cap \{\wt{Y}_2, \dots, \wt{Y}_r\} = \emptyset$. Each element of $\setrow Yn$ is an effective Cartier divisor so again we can find a cover of $X$ by open affines $\{U_\alpha = \Spec R_\alpha\}$ such that for each $\alpha$, and each $Y \in \setrow Yn$, we have a non-zerodivisor $y \in R_\alpha$ which cuts out $Y$ in $U_\alpha$.  Choose such a $U = \Spec R$ and let $y_i \in R$ be a local equation for $Y_i$.  We can then write the blowup of this chart
	 		\[
		 		\wt U = \Proj R[a_1, \dots, a_r]/(\{ a_j y_i - a_i y_j \;|\; i \neq j \})
		 	\] 
		 	(again see \cite[B.6.10]{Ful98}). As before, we work affine-locally in 
		 	\[
			 	D(a_1) = \Spec R[a_2, \dots, a_r]/(\{ y_i - a_i y_1 \;|\; i \neq 1 \})
			\]
			where we have that $E$ is cut out by $y_1$, and $\wt Y_i$ is cut out by $a_i$ for each $i >2$.

			We can use induction on $r$ with the base case \autoref{blowupProp}. Then assume $\{E, \wt{Y}_1, \dots, \wt{Y}_n\}$ has c.i.~crossings in the blowup
			$\beta': X' \rightarrow X$ centered at $\cap_{i=1}^{r-1} Y_i$.
			If we let $U'$ denote the blowup of the chart $U$ by $\beta'$ then
			\[
				U' = \Proj R[\row a{r-1}]/(\{ a_j y_i - a_i y_j \;|\; i \neq j \})
			\] 
			with corresponding open chart
			\[
				D'(a_1) = \Spec R[\row a{r-1}]/(\{ y_i - a_i y_1 \;|\; i = 2,\dots, r-1 \}).
			\]		
			Now we observe that $D(a_1)$ is related to $D'(a_1)$ by
			\begin{align*}
				D(a_1) \cap \mathbb{V}(a_r) &= \Spec (R[a_2, \dots, a_r]/(\{ y_i - a_i y_1 \;|\; i \neq 1 \})) /(a_r) \\
				&= \Spec R[a_2, \dots, a_{r-1}] / (\{ y_i - a_i y_1 \;|\; i = 2, \dots, r-1 \} \cup \{y_r\}) \\
				&= D'(a_1) \cap \mathbb{V}(y_r).
			\end{align*}
			Recall that $\wt Y_r \notin A$. The inductive hypothesis says that we can find a regular sequence $\row sk$ corresponding to $A$, where $k = \#A$. Then as usual we can localize at a point and rearrange to get that $y_r, \row sk$ is a regular sequence at each point $p \in D'(a_1)$.  Then $\row sk$ is a regular sequence at each point in $D'(a_1)/(y_r) = D(a_1)/(a_r)$, so $a_r, \row sk$ is a regular sequence at each point of $D(a_1)$.
	 	\end{proof}
	\end{prop}

\section{Proof of main theorem}\label{sec:mainThmProof}

	In \cite{Gow05}, Goward defines an invariant $(\sigma, \tau)$ on the divisors 
	in question.  We adopt those definitions for the new context here:

	\begin{defn}\label{sigmaDefn}
		Let $D_1 = \sum_{i=1}^n a_i Y_i$ and $D_2 = \sum_{i=1}^n b_i Y_i$ where $\setrow Yn$ has c.i.~crossings in $X$ and $a_i,b_j \geq 0$. We define
		\[ 
			\sigma_{ij}(D_1,D_2) = \begin{cases}
				\max \{ (|a_i - b_i|,|a_j-b_j|), (|a_j-b_j|,|a_i-b_i|) \} & \parbox[t]{.25\textwidth}{\footnotesize if $a_i-b_i$ and $a_j-b_j$\\ have opposite signs,}\\
				(-\infty, -\infty) & \text{\small otherwise,}
			\end{cases}
		\]
		where the $\max$ is taken lexicographically.
		Now we can define
		\begin{align*}
			\sigma(D_1,D_2) &= \max \{ \sigma_{ij}(D_1,D_2) \;|\; Y_i \cap Y_j \neq 0,\; i \neq j \} \\
			\tau(D_1,D_2) &= \# \{ (i,j) \;|\; \sigma_{ij}(D_1,D_2) = \sigma(D_1,D_2),\; i \leq j \}
		\end{align*}
		so $\sigma(D_1,D_2)$ takes the value of the worst intersection in the support, and $\tau(D_1,D_2)$ counts how many intersections share this value.
	\end{defn}

	These invariants are calculated for divisors on $X$ and, after blowing up, for their pullbacks.  We show that these calculations go the same way as in the simple normal crossings context and then outline the steps of the proof.

	\begin{prop}\label{sigmaDetectsPrincipalityProp}
		Let $D_1, D_2$ be as defined above.
		Then $\scrI_{D_1} + \scrI_{D_2}$ is principal at $p \in X$ if and only if 
		$\sigma_{ij}(D_1,D_2) = (-\infty,-\infty)$ whenever $p \in Y_i \cap Y_j$.

		\begin{proof}
			Let $\row xn $ be a regular sequence at $p \in X$ corresponing to $\row Yn$. Then let $f_1 = u_1 x_1^{a_1} \dots x_n^{a_n}$ and $f_2 = u_2 x_1^{b_1} \dots x_n^{b_n}$ be local equations for $D_1$ and $D_2$, where $u_i \in \calO_{X,p}$ are units. Suppose $\scrI_{D_1} + \scrI_{D_2}$ is not principal at $p$.  Then $(f_1,f_2)$ is not principal, so we have some $i,j$ such that $a_i - b_i$ and $a_j - b_j$ have opposite signs. Thus $\sigma_{ij}(D_1,D_2) > 0$. On the other hand, suppose $p \in Y_i \cap Y_j$ and $\sigma_{ij}(D_1,D_2) > 0$. Then $a_i - b_i$ and $a_j - b_j$ have opposite signs, so $(f_1,f_2)$ is not principal and thus $\scrI_{D_1} + \scrI_{D_2}$ is not principal at $p$.
		\end{proof}
	\end{prop}

	As a result of this proof, we see that if $\scrI_{D_1} + \scrI_{D_2}$ is principal at $p$, then $(\scrI_{D_1} + \scrI_{D_2})_p = (\scrI_{D_i})_p$ for some $i \in \{1,2\}$. So by induction, we have that  if $\scrI_{D_1} + \dots + \scrI_{D_h}$ is principal at $p$, then $(\scrI_{D_1} + \dots + \scrI_{D_h})_p = (\scrI_{D_i})_p$ for some $i \in \{1, \dots, h\}$.

	We have left to show that blowing up at the chosen codimension 2 centers strictly reduces the invariant $(\sigma,\tau)$ and then that such blowups can be taken successively until $(\sigma,\tau) = (-\infty, -\infty)$.

	\begin{prop}\label{sigmaDecreasesProp}
		Let $D_1, D_2$ be as defined above. Suppose we have $(i,j)$ 
		such that $\sigma_{ij}(D_1,D_2) = \sigma(D_1,D_2) > (-\infty, -\infty)$ and let $\beta: \wt X \rightarrow X$ be the blowup of $X$ centered at $Y_i \cap Y_j$.
		Then 
		\[
			(\sigma(D_1,D_2), \tau(D_1,D_2)) 
			> (\sigma(\beta^*D_1,\beta^*D_2),\tau(\beta^*D_1,\beta^*D_2).
		\]

		\begin{proof}[Sketch of proof]
			Assume that $(i,j) = (1,2)$ so that the blowup is centered at $Y_1 \cap Y_2$. The proof relies on calculations of $\sigma_{ij}(\beta^* D_1, \beta^* D_2)$, which depend only on the coefficients $a_i, b_i, a_j, b_j$.  The details can be found in the proof of \cite[Thm 1]{Gow05}.  Here we verify only that the calculations of $\sigma_{ij}(\beta^* D_1, \beta^* D_2)$ from the s.n.c.~case still go through with c.i.~crossings.

			Let $E \subset \wt X$ denote the exceptional divisor.  Then
			\[ 
				\beta^*Y_i = \begin{cases}
					\wt Y_i + E & \text{if}~ i = 1,2 \\
					\wt Y_i & \text{if}~ i > 2
				\end{cases}
			\]
			as seen from the work done in \autoref{blowupBaseCaseProp}.  Thus
			\[ 
				\beta^*{D_1} = (a_1 + a_2) E + \sum_i a_i \wt Y_i
			\]
			and
			\[
				\beta^*{D_2} = (b_1 + b_2) E + \sum_i b_i \wt Y_i. \qedhere
			\]
		\end{proof}
	\end{prop}

	\newtheorem*{mainThm}{Theorem \ref{mainThm}}
	\begin{mainThm}[Main Theorem]
		\stateMainThm
		\begin{proof}
			By the previous remarks, we can assume $h = 2$.  Let $\row Yn$ be the 
			divisors in the support of $D_1 + D_2$.  If 
			$\sigma (D_1,D_2) = (-\infty, -\infty)$ we are done, so assume not. 
			Then let 
			\[
				(i,j) = 
				\max \{ (k,l) \;|\; \sigma_{kl}(D_1, D_2) = \sigma(D_1,D_2) \}
			\]
			where the $\max$ is taken lexicographically. If we take the blowup 
			$\beta_1: X_1 \rightarrow X$ centered at $Y_i \cap Y_j$, then 
			\autoref{sigmaDecreasesProp} gives that
			\[
				(\sigma(D_1,D_2), \tau(D_1,D_2)) 
				> (\sigma(\beta^*D_1,\beta^*D_2),\tau(\beta^*D_1,\beta^*D_2).
			\]
			Now we note that $\beta_1^* D_1$, $\beta_1^* D_2$ are divisors supported on $\{E, \row {\wt Y}n \}$ and have c.i.~crossings by \autoref{blowupProp}.  Thus, $(\scrI_{D_1} + \scrI_{D_2}) \calO_{\wt X}$ defines a c.i.~monomial subscheme of $\wt X$.

			We can repeat this process, with $(\sigma,\tau)$ decreasing at each iteration. Since $(\sigma, \tau)$ takes values in $\bbN^2 \times \bbN$, we must get $(\sigma, \tau) = (-\infty, -\infty)$ after finitely many steps, yielding the desired sequence of blowups.
		\end{proof}
	\end{mainThm}

% \section{Conclusions}\label{sec:conclusions}

%%%%%%%%%%%%%%%%%%%%%%%%%%
% Bibliography
%%%%%%%%%%%%%%%%%%%%%%%%%%
\bibliography{sp_refs}{}
\bibliographystyle{amsalphaeprint}

% \printbibliography
%%%%%%%%%%%%%%%%%%%%%%%%%%

\end{document}